\documentclass[preprint,12pt]{elsarticle}

\usepackage{amssymb}
\usepackage{latexsym}
\usepackage{amscd}
\usepackage[centertags]{amsmath}

\journal{arXiv}

\newtheorem{thm}{Theorem}

\newcommand{\Enum}{\mathbb{E}}

\begin{document}

\begin{frontmatter}



\title{Range-Renewal Structure of I.I.D. Samplings}


\author{Li-Hong Wang and Jian-Sheng Xie}
\address{School of Mathematical Sciences, Fudan University,
Shanghai 200433, China}
\address{E-mail: jiansheng.xie@gmail.com}

\begin{abstract}
In this note the range-renewal structure of general independent and identically distributed 
samplings is discussed, which is a natural extension of the result in Chen et al (arXiv:1305.1829).
\end{abstract}

\end{frontmatter}


%

Let $\xi:=\{\xi_n: n \geq 1\}$ be an i.i.d. symbol sequence (with common distribution law
$\pi$) and let $R_n$ be the number of distinct elements among the first $n$ elements of the
process $\xi$. For latter use, put
\begin{equation}
N_n (x) :=\sum_{k=1}^n 1_{\{\xi_k=x\}},
\end{equation}
which is the number of visiting times (visiting intensity) at $x$ of the random sequence up
to time $n$. Then put for each $\ell \geq 1$
\begin{equation}
R_{n, \, \ell} :=\sum_x 1_{\{N_n (x)=\ell\}}, \quad R_{n, \, \ell+} :=\sum_x
1_{\{N_n (x) \geq \ell\}};
\end{equation}
Then $R_{n, \, \ell}$ (resp. $R_{n, \, \ell+}$) is the number of distinct states which have
been visited at
exactly (resp. at least) $\ell$ times in the first $n$ steps. Obviously
$$
R_n=R_{n,\, 1+}=\sum_{\ell=1}^n R_{n, \, \ell}, \quad R_{n, \,
\ell+}=\sum_{k=\ell}^n R_{n, \, k}.
$$
Chen et al (arXiv:1305.1829) discovered a range-renewal structure for the i.i.d. samplings of
\textit{regular} discrete distributions as the following. For general discrete distribution
$\pi$ which has infinite many atoms, one has
\begin{equation}
\lim_{n \to \infty} \frac{R_n}{\Enum R_n}=1
\end{equation}
almost surely. When the distribution $\pi$ is in fact regular (therefore with an index
$\gamma=\gamma (\pi) \in [0, 1]$), it is proved that for all $k \geq 1$
\begin{eqnarray}
\label{eq: non-critical} \lim_{n \to \infty} \frac{R_{n, \, k}}{R_n} &=r_k (\gamma), \; 
\lim \limits_{n \to \infty} \frac{R_{n, \, k}}{R_{n, \, k+}} &=\frac{\gamma}{k} 
\quad \hbox{ if } 0< \gamma <1\\
\label{eq: sub-critical} \lim_{n \to \infty} \frac{R_{n, \, k}}{R_n} &=0,  \qquad  \lim 
\limits_{n \to \infty} \frac{R_{n, \, k}}{R_{n, \, k+}} &=0 \quad \hbox{ if } \gamma=0\\
\label{eq: sup-critical} \lim_{n \to \infty} \frac{R_{n, \, 1}}{R_n} &=1, \quad \lim 
\limits_{n \to \infty} \frac{R_{n, \, k+1}}{R_{n, \, 2+}} &=\frac{1}{k (k+1)} 
\hbox{ if } \gamma=1,
\end{eqnarray}
where in Eq. (\ref{eq: non-critical}) $r_k (\gamma) := 
\frac{\gamma \cdot \Gamma (k-\gamma)}{k! \cdot \Gamma (1-\gamma)}$.

In this note we would like to discuss what happens for general distributions which are not
purely discrete. Let $\pi$ be such a distribution on an abstract phase space $E$. Let
$$
\Delta :=\{x \in E: \pi (\{x\})>0\}
$$
be the set of atoms of $\pi$. If $\Delta =\emptyset$, it's clear that $R_n=R_{n, 1}=n$
almost surely. Therefore let's first assume that $\pi$ has infinitely many atoms and that
the conditional distribution $\pi^* :=\frac{\pi|_{\Delta}}{\pi (\Delta)}$ is regular in the
sense explained in \cite{CXY} with an index $\gamma^*:=\gamma (\pi^*)$. We would prove the
following result.
\begin{thm}\label{thm: main}
Let $\pi$ be probability distribution on $E$ with atoms set $\Delta$ of infinite cardinality.
Suppose that the conditional distribution $\pi^* :=\frac{\pi|_{\Delta}}{\pi (\Delta)}$ is
regular in the sense explained in \cite{CXY} with an index $\gamma^*:=\gamma (\pi^*)$. Then
for a sequence of i.i.d. sampling $\xi=\{\xi_n: n \geq 1\}$ with common distribution $\pi$,
we have
\begin{equation}\label{eq: 1}
\lim_{n \to +\infty} \frac{R_n}{n} = \lim_{n \to +\infty} \frac{R_{n, \,1}}{n}
= 1-\pi (\Delta)
\end{equation}
and for any $k \geq 2$
\begin{equation}\label{eq: 2}
\lim_{n \to +\infty} \frac{R_{n, \, k}}{R_{n, \, 2+}} =\left\{
\begin{array}{rcl}
\frac{r_k (\gamma^*)}{1-\gamma^*}, \qquad &\hbox{ if }& \gamma^* \in (0, 1)\\
0, \qquad  &\hbox{ if }& \gamma^*=0\\
\frac{1}{k (k-1)}, \qquad  &\hbox{ if }& \gamma^*=1.
\end{array}
\right.
\end{equation}

\begin{equation}\label{eq: 3}
\lim_{n \to +\infty} \frac{R_{n, \, k}}{R_{n, \, k+}} =\left\{
\begin{array}{rcl}
\frac{\gamma^*}{k}, \qquad &\hbox{ if }& \gamma^* \in (0, 1]\\
0, \qquad  &\hbox{ if }& \gamma^*=0.
\end{array}
\right.
\end{equation}
almost surely.
\end{thm}

The proof of the above result is simple. Let's first define two sequence of stopping times 
as the following.
\begin{eqnarray*}
\tau_0 :=0, & \tau_k :=\inf\{n >\tau_{k-1}: \xi_n \in \Delta\}, \quad k =1, 2, \cdots,\\
s_0 :=0, & s_k :=\inf\{n >s_{k-1}: \xi_n \in \Delta^c\}, \quad k =1, 2, \cdots.
\end{eqnarray*}
Then for any $k \geq 1$ we put
$$
X_k :=\xi_{\tau_k}, \quad Y_k :=\xi_{s_k}, \quad Z_k :=1_{\{\xi_k \in \Delta^c\}}.
$$
It's easy to see that each of the processes $X=\{X_n: n \geq 1\}, Y=\{Y_n: n \geq 1\}, 
Z=\{Z_n: n \geq 1\}$ is i.i.d. samplings with corresponding common distributions being $\pi^* 
:=\frac{\pi|_{\Delta}}{\pi (\Delta)}, \frac{\pi|_{\Delta^c}}{\pi (\Delta^c)}$ and the classical 
Bernoulli distribution $B (1, p)$ which assigns probability $p=\pi (\Delta^c)$ to $1$ and 
probability $q=1-p$ to $0$. 
Let 
$$
N_n (Z) :=\sum_{k=1}^n Z_k.
$$
Let $R_n=R_n (\xi)$ be the number of distinct states among $\xi_1, \cdots, \xi_n$; we emphasis 
here the obvious dependence of the number $R_n$ on the process $\xi$. Other notations 
have similar meanings. Then it is easy to see that
\begin{eqnarray*}
R_n (\xi) &=& N_n (Z) +R_{n-N_n (Z)} (X), \\
R_{n, \, 1} (\xi) &=& N_n (Z)+R_{n-N_n (Z), \, 1} (X)\\
R_{n, \, k} (\xi) &=& R_{n-N_n (Z), \, k} (X), \; k=2, 3, \cdots
\end{eqnarray*}
almost surely. Since we have a strong law of large numbers here
$$
\lim_{n \to \infty} \frac{N_n (Z)}{n}=p=\pi (\Delta^c)=1-\pi (\Delta),
$$
the results in Theorem \ref{thm: main} follows easily from eq. (\ref{eq: non-critical})--(\ref{eq: sup-critical}).

We remark here that, when the atom set $\Delta$ of the distribution $\pi$ is non-empty and
finite, the result stated in Theorem \ref{thm: main} still holds true with $\gamma^* :=0$; 
The proof is almost the same with necessary modifications.

\noindent{\sl \textbf{Acknowledgements} \quad} This work is in part supported by NSFC
No. 11271077 and the Laboratory of Mathematics for Nonlinear Science, Fudan University.





\bibliographystyle{elsarticle-num}



\end{document}